\documentclass[preprint,3p,times]{elsarticle}

\usepackage{amssymb}
\usepackage{amsmath}
\usepackage[all,cmtip]{xy}
\usepackage{amsthm}
\usepackage{color}
\usepackage{enumitem}
\usepackage{mathtools}
\usepackage{tikz}
\usepackage{amsmath, mathrsfs}
\usepackage{appendix}
\usepackage[colorlinks=true]{hyperref}
\usepackage[hyphenbreaks]{breakurl}
\input diagxy

\theoremstyle{plain}
 \newtheorem{thm}{Theorem}[section]
 \newtheorem{lem}[thm]{Lemma}
 \newtheorem{prop}[thm]{Proposition}
 \newtheorem{cor}[thm]{Corollary}
\theoremstyle{definition}
 
 \newtheorem{exmp}[thm]{Example}
 \newtheorem{rem}[thm]{Remark}

\DeclareMathAlphabet{\mathcal}{OMS}{cmsy}{m}{n}

\DeclareMathOperator{\ev}{ev}

\DeclareMathOperator{\ob}{ob}

\makeatletter
\def\ps@pprintTitle{%
\let\@oddhead\@empty
\let\@evenhead\@empty
\def\@oddfoot{\centerline{\thepage}}%
\let\@evenfoot\@oddfoot}
\makeatother

\def\oto{{\bfig\morphism<180,0>[\mkern-4mu`\mkern-4mu;]\place(86,0)[\circ]\efig}}

\newcommand{\ra}{\rightarrow}

\newcommand{\bv}{\bigvee}

\newcommand{\dv}{\dashv}

\newcommand{\nat}{\natural}

\newcommand{\dbv}{\displaystyle\bv}

\renewcommand{\phi}{\varphi}
\newcommand{\al}{\alpha}
\newcommand{\be}{\beta}

\newcommand{\ga}{\gamma}

\newcommand{\lam}{\lambda}

\newcommand{\Om}{\Omega}

\newcommand{\CC}{\mathcal{C}}

\newcommand{\sC}{\mathsf{C}}
\newcommand{\sD}{\mathsf{D}}
\newcommand{\sE}{\mathsf{E}}

\newcommand{\sM}{\mathsf{M}}

\newcommand{\sQ}{\mathsf{Q}}
\newcommand{\sR}{\mathsf{R}}

\newcommand{\sY}{\mathsf{Y}}

\newcommand{\bbN}{\mathbb{N}}

\newcommand{\Fr}{\mathfrak{r}}

\newcommand{\Fy}{\mathfrak{y}}

\newcommand{\Set}{\mathbf{Set}}
\newcommand{\CcSet}{\mathbf{CcSet}}

\newcommand{\QSet}{\sQ\text{-}\Set}

\newcommand{\RSet}{[0,1]_*\text{-}\Set}

\newcommand{\RSetM}{[0,1]_*\text{-}\underline{\Set}}
\newcommand{\RCcSetM}{[0,1]_*\text{-}\underline{\CcSet}}
\newcommand{\RLCcSetM}{[0,1]_{*_{\L}}\text{-}\underline{\CcSet}}

\newcommand{\sCd}{\sC^{\dag}}

\newcommand{\Fyd}{\Fy^{\dag}}

\newcommand{\CdX}{\sCd X}

\newcommand{\DQ}{\sD\sQ}

\renewcommand{\leq}{\leqslant}

\numberwithin{equation}{section}

\allowdisplaybreaks

\begin{document}

\begin{frontmatter}



\title{Cartesian closedness of the category of real-valued sets, I}


\author{Lili Shen}
\ead{shenlili@scu.edu.cn}

\author{Jian Zhang\corref{cor}}
\ead{zhangjian.2025@qq.com}

\cortext[cor]{Corresponding author.}
\address{School of Mathematics, Sichuan University, Chengdu 610064, China}

\begin{abstract}
Let $[0,1]_*$ be the unit interval $[0,1]$ equipped with a continuous t-norm $*$.  It is shown that the category of $[0,1]_*$-sets is cartesian closed if, and only if, $*$ is the minimum t-norm on $[0,1]$.
\end{abstract}


\end{frontmatter}




\section{Introduction}

Building on the theory of frame-valued sets of Higgs \cite{Higgs1970,Higgs1984} and Fourman–Scott \cite{Fourman1979}, the theory of quantale-valued sets developed by H{\"o}hle and his collaborators \cite{Hoehle1991,Hoehle1992,Hoehle1995b,Hoehle1998,Hoehle2005,Hoehle2011a} has been influential in the categorical foundations of fuzzy sets \cite{Zadeh1965}. Given a frame $\Om$, an \emph{$\Om$-set} \cite{Fourman1979,Borceux1994c} consists of a (crisp) set $X$ and a map
\[\al\colon X\times X\to\Om\]
such that
\begin{equation} \label{Om-set-def}
\al(x,y)=\al(y,x)\quad\text{and}\quad\al(y,z)\wedge\al(x,y)\leq\al(x,y)
\end{equation}
for all $x,y,z\in X$, where $\al(x,y)$ measures the degree of $x$ being equal to $y$, and \eqref{Om-set-def} refers to the symmetry and transitivity of the $\Om$-valued equality $\al$. It is well known that the category
\[\Om\text{-}\Set\]
of $\Om$-sets and their morphisms is a topos (see \cite[Theorem 5.9 and Proposition 9.2]{Fourman1979}). Consequently, $\Om\text{-}\Set$ enjoys many desirable properties, including cartesian closedness and the existence of a subobject classifier \cite{Borceux1994c}. 

From the viewpoint of \emph{enriched category theory} \cite{Kelly1982}, every frame $\Om$ gives rise to a bicategory $\sD\Om$ \cite{Walters1981}, and $\sD\Om$ is actually a \emph{quantaloid} \cite{Rosenthal1996,Stubbe2005,Stubbe2014}. $\Om$-sets are exactly symmetric categories enriched in the quantaloid $\sD\Om$, and morphisms between $\Om$-sets are left adjoint $\sD\Om$-distributors.

If we consider a unital and involutive \emph{quantale} \cite{Rosenthal1990,Mulvey1992} $\sQ$ as the table of truth values, the theory of $\sQ$-sets can be established in the same way. Explicitly, we may construct a quantaloid $\DQ$ \cite{Hoehle2011a,Pu2012,Stubbe2014}, and define $\sQ$-sets as symmetric $\DQ$-categories \cite{Hoehle2011a,Pu2012} (cf. Remark \ref{D[0,1]-rem}). However, the category 
\[\QSet\]
of $\sQ$-sets and their morphisms need not be a topos. More precisely, in the case that $\sQ$ is a commutative, unital and divisible, Hu--Shen \cite{Hu2025} proved that $\QSet$ is a topos if, and only if, $\sQ$ is a frame.

Although $\QSet$ generally fails to be a topos, it is reasonable to ask which topos-like features it may still possess. In this paper we initiate a study of the cartesian closedness of $\QSet$. Given the difficulty in treating an arbitrary quantale $\sQ$, we begin with the special case 
\[\sQ=[0,1]_*,\] 
where $*$ is a \emph{continuous t-norm} \cite{Klement2000,Klement2004b,Alsina2006} on the unit interval $[0,1]$. Using the linearity of $[0,1]$ and the continuity of $*$, we derive our main result (see Theorem \ref{QSet-topos-frame}):
\begin{itemize}
\item The category $\RSet$ is cartesian closed if, and only if, $*$ is the minimum t-norm on $[0,1]$.
\end{itemize}
Therefore, $\RSet$ also fails to be cartesian closed when $[0,1]_*$ is not a frame.

To make this paper accessible to readers who may not be familiar with enriched category theory, we present the theory of $[0,1]_*$-sets in fairly elementary terms, even though it can be formulated directly via quantaloid-enriched categories (cf. Remark \ref{D[0,1]-rem}). Moreover, it is tempting to ask whether the methods adopted in this paper can extended to an arbitrary quantale $\sQ$. Unfortunately, this is not possible due to the lack of linearity and continuity in the general setting. In a subsequent work, we will investigate the cartesian closedness of the category $\RSet$ while weakening the assumption on $*$ to left-continuity.

\section{Real-valued sets}

Throughout this paper, we consider $[0,1]_*$ as the table of truth-values, where $*$ is a continuous t-norm on the unit interval $[0,1]$. More generally, a binary operation $*$ on an interval $[a,b]$ is a \emph{continuous t-norm} \cite{Klement2000,Klement2004b,Alsina2006}, denoted by $[a,b]_*$, if
\begin{itemize}
\item $([a,b],*,b)$ is a commutative monoid,
\item $p*q\leq p'*q'$ if $p\leq p'$ and $q\leq q'$ in $[a,b]$, and
\item $*\colon [a,b]\times[a,b]\to[a,b]$ is a continuous function (with respect to the usual topology).
\end{itemize}


Note that for each $p\in[a,b]$, there is a Galois connection
\[(p*-)\dv(p\ra -)\colon[a,b]\to[a,b]\]
satisfying
\[p*q\leq r\iff q\leq p\ra r\]
for all $p,q,r\in[a,b]$. 

Let $p\in[a,b]$. We say that $p$ is \emph{idempotent} if $p*p=p$. In this case, it is easy to see that
\begin{equation} \label{p-idem-q}
p*q=p\wedge q\quad\text{and}\quad p\ra r=r
\end{equation}
for all $q\in[a,b]$, $r\in[a,p)$. 


\begin{exmp} \label{MLP-def}
The following continuous t-norms on the unit interval $[0,1]$ are the most prominent ones:
\begin{enumerate}[label=(\arabic*)]
\item \label{MLP-def:M} The \emph{minimum t-norm} $[0,1]_{\wedge}$ is given by the meet of real numbers, in which every $q\in(0,1)$ is idempotent, and
\[p\ra q=\begin{cases}
1 & \text{if}\ p\leq q,\\
q & \text{if}\ p>q.
\end{cases}\]
\item \label{MLP-def:P} The \emph{product t-norm} $[0,1]_{\times}$ is given by the usual multiplication of real numbers, in which every $q\in(0,1)$ is non-idempotent, and
\[p\ra q=1\wedge\dfrac{q}{p}.\]
\item \label{MLP-def:L} The \emph{{\L}ukasiewicz t-norm} $[0,1]_{*_{\L}}$ is given by $p*_{\text{\L}} q=0\vee(p+q-1)$ for all $p,q\in[0,1]$, in which every $q\in(0,1)$ is non-idempotent, and
\[p\ra q=1\wedge(1-p+q).\]
\end{enumerate}
\end{exmp}

In fact,
\[([a,b],*,b)\]
is a commutative, unital and divisible \emph{quantale} \cite{Rosenthal1990}. 

\begin{prop} \label{divisible} (See \cite[Proposition 2.1]{Pu2012}.)
The following identities hold in every continuous t-norm $[a,b]_*$:
\begin{enumerate}[label={\rm(\arabic*)}]
\item \label{divisible:u-q} $u=q*(q\ra u)$ whenever $u\leq q$ in $[a,b]$.
\item \label{divisible:u-v-q} $v*(q\ra u)=(q\ra v)*u$ whenever $u,v\leq q$ in $[a,b]$.
\item \label{divisible:p-q} $p\wedge q=p*(p\ra q)$ for all $p,q\in[a,b]$.
\end{enumerate}
\end{prop}

We say that continuous t-norms $[a_1,b_1]_*$ and $[a_2,b_2]_{\bullet}$ are \emph{isomorphic} if the quantales $([a_1,b_1],*,b_1)$ and $([a_2,b_2],\bullet,b_2)$ are isomorphic; that is, if there exists an order isomorphism
\[f\colon[a_1,b_1]\to[a_2,b_2]\] 
such that 
\[f\colon([a_1,b_1],*,b_1)\to([a_2,b_2],\bullet,b_2)\] 
is an isomorphism of monoids.

It is well known \cite{Faucett1955,Mostert1957,Klement2000,Klement2004b,Alsina2006} that every continuous t-norm $*$ on $[0,1]$ can be written as an \emph{ordinal sum} of the minimum, the product, and the {\L}ukasiewicz t-norm. Explicitly:

\begin{lem} \label{t-norm-rep} {\rm\cite{Klement2000,Klement2004b,Alsina2006}}
For each continuous t-norm $[0,1]_*$, the set of non-idempotent elements of $*$ in $[0,1]$ is a union of countably many pairwise disjoint open intervals 
\[\{(p_i,q_i)\mid 0<p_i<q_i<1,\ i\in I,\ I\ \text{is countable}\},\]
and for each $i\in I$, the continuous t-norm $[p_i,q_i]_*$ obtained by restricting $*$ to $[p_i,q_i]$ is either isomorphic to the product t-norm $[0,1]_{\times}$ or isomorphic to the {\L}ukasiewicz t-norm $[0,1]_{*_{\L}}$. 
\end{lem}

Following the definition of \emph{quantale-valued set} \cite{Hoehle1992,Hoehle1995b,Hoehle2011a,Pu2012}, by a \emph{$[0,1]_*$-set} (or, \emph{$[0,1]_*$-valued set}) we mean a (crisp) set $X$ equipped with a map
\[\al\colon X\times X\to[0,1],\]
such that
\begin{enumerate}[label=(S\arabic*)]
\item \label{S1} $\al(x,y)\leq\al(x,x)\wedge\al(y,y)$,
\item \label{S2} $\al(x,y)=\al(y,x)$,
\item \label{S3} $\al(y,z)*(\al(y,y)\ra\al(x,y))\leq\al(x,z)$
\end{enumerate}
for all $x,y,z\in X$. 

\begin{rem}
A $[0,1]_*$-set $(X,\al)$ may be viewed as a set $X$ equipped with a \emph{$[0,1]_*$-valued equality} (or, \emph{$[0,1]_*$-valued similarity}) $\al$ \cite{Hoehle2011a,Lai2020}. The value $\al(x,y)$ is interpreted as the extent of $x$ being equal to $y$, and $\al(x,x)$ represents the extent of existence of $x$ (since every entity is supposed to be equal to itself). Therefore:
\begin{itemize}
\item \ref{S1} says that $x$ is equal to $y$ only if both $x$ and $y$ exist; that is, \emph{equality} implies \emph{existence}.
\item \ref{S2} says that if $x$ is equal to $y$, then $y$ is equal to $x$.
\item \ref{S3} says that if $x$ is equal to $y$, and there exists $y$ such that $y$ is equal to $z$, then $x$ is equal to $z$. 
\end{itemize}
\end{rem}

Let $(X,\al)$ be a $[0,1]_*$-set. For $x,y\in X$, we write $x\cong y$ if
\begin{equation} \label{x-cong-y-def}
\al(x,x)=\al(y,y)=\al(x,y),
\end{equation}
and we say that $(X,\al)$ is \emph{separated} if 
\[x\cong y\iff x=y.\]
Moreover, we denote by
\begin{equation} \label{X_q-def}
X_q=\{x\in X\mid\al(x,x)=q\}
\end{equation}
the slice of $X$ consisting of elements whose extent of existence is $q$.

A \emph{morphism} 
\[\phi\colon(X,\al)\oto(Y,\be)\] 
of $[0,1]_*$-sets is a function
\[\phi\colon X\times Y\to[0,1],\]
such that
\begin{enumerate}[label=(M\arabic*)]
\item \label{M1} $\phi(x,y)\leq\al(x,x)\wedge\be(y,y)$,
\item \label{M2} $(\be(y,y)\ra\be(y,y'))*\phi(x,y)*(\al(x,x)\ra\al(x',x))\leq\phi(x',y')$,
\item \label{M3} $\phi(x,y')*(\al(x,x)\ra\phi(x,y))\leq\be(y,y')$,
\item \label{M4} $\al(x,x')\leq\dbv\limits_{y\in Y}\phi(x',y)*(\be(y,y)\ra\phi(x,y))$
\end{enumerate}
for all $x,x'\in X$, $y,y'\in Y$; its composite with another morphism $\psi\colon(Y,\be)\oto(Z,\ga)$ is given by
\begin{equation} \label{psi-circ-phi-def}
\psi\circ\phi\colon(X,\al)\oto(Z,\ga),\quad (\psi\circ\phi)(x,z)=\bv_{y\in Y}\psi(y,z)*(\be(y,y)\ra\phi(x,y)),
\end{equation}
with
\begin{equation} \label{al-id}
\al\colon(X,\al)\oto(X,\al)
\end{equation}
playing as the identity morphism for the composition. The category of $[0,1]_*$-sets and their morphisms is denoted by
\[\RSet.\]

\begin{rem} \label{OmSet}
If $*$ is the minimum t-norm $\wedge$ on $[0,1]$ (see Example \ref{MLP-def}\ref{MLP-def:M}), then the category
\[[0,1]_{\wedge}\text{-}\Set\]
is a special case of the the well-known category $\Om\text{-}\Set$, where $\Om$ is a frame; see \cite{Fourman1979} and \cite[Sections 2.8 and 2.9]{Borceux1994c}. In this case, $[0,1]_{\wedge}\text{-}\Set$ is equivalent to the category $\mathbf{Sh}([0,1]_{\wedge})$ of sheaves on $[0,1]_{\wedge}$ (cf. \cite[Theorem 5.9]{Fourman1979} and \cite[Theorem 2.9.8]{Borceux1994c}). Therefore, $[0,1]_{\wedge}\text{-}\Set$ is a topos (see \cite[Proposition 9.2]{Fourman1979} and \cite[Example 5.2.3]{Borceux1994c}) and, in particular, it is a cartesian closed category.
\end{rem}

We may also consider \emph{monotone functions} 
\[f\colon(X,\al)\to(Y,\be)\] 
between $[0,1]_*$-sets, which is a function $f\colon X\to Y$ such that
\begin{equation} \label{monotone-function-def}
\al(x,x)=\be(fx,fx)\quad\text{and}\quad\al(x,x')\leq\be(fx,fx')
\end{equation}
for all $x,x'\in X$. The category of $[0,1]_*$-sets and monotone functions is denoted by
\[\RSetM.\]

\begin{rem} \label{D[0,1]-rem}
It is noteworthy to point out that every continuous t-norm $[0,1]_*$ gives rise to a \emph{quantaloid} $\sD[0,1]_*$ \cite{Hoehle2011,Pu2012,Stubbe2014,Lai2020} consisting of the following data:
\begin{itemize}
\item Objects of $\sD[0,1]_*$ are elements of $[0,1]$.
\item For $p,q\in[0,1]$, the hom-set $\sD[0,1]_*(p,q)=\{u\in[0,1]\mid u\leq p\wedge q\}$.
\item The composite of $u\in\sD[0,1]_*(p,q)$ and $v\in\sD[0,1]_*(q,r)$ is given by  
      \[v\circ u=v*(q\ra u).\]
\item The identity $\sD[0,1]_*$-arrow of $\sD[0,1]_*(q,q)$ is $q$ itself.
\item Each hom-set $\sD[0,1]_*(p,q)$ is equipped with the order inherited from $[0,1]$.
\end{itemize}
From the viewpoint of enriched category theory, the foregoing notions of $[0,1]_*$-sets can be formulated within the framework of $\sD[0,1]_*$-categories:
\begin{itemize}
\item A $[0,1]_*$-set is a symmetric \emph{$\sD[0,1]_*$-category}, where the type function is given by $x\mapsto\al(x,x)$.
\item A morphism of $\sD[0,1]_*$-sets is a \emph{left adjoint $\sD[0,1]_*$-distributor} between $\sD[0,1]_*$-categories;
\item A monotone function between $\sD[0,1]_*$-sets is a \emph{$\sD[0,1]_*$-functor} between $\sD[0,1]_*$-categories.
\end{itemize}
Therefore:
\begin{itemize}
\item $\RSet$ is the category of symmetric $\sD[0,1]_*$-categories and left adjoint $\sD[0,1]_*$-distributors.
\item $\RSetM$ is the category of symmetric $\sD[0,1]_*$-categories and $\sD[0,1]_*$-functors.
\end{itemize}
\end{rem}

Every monotone function $f\colon(X,\al)\to(Y,\be)$ induces a morphism 
\begin{equation} \label{f-graph-def}
f_{\nat}\colon(X,\al)\oto(Y,\be),\quad f_{\nat}(x,y)=\be(fx,y)
\end{equation}
of $[0,1]_*$-sets, called the \emph{graph} of $f$, whose opposite
\[f^{\nat}\colon(Y,\be)\oto(X,\al),\quad f^{\nat}(y,x)=\be(y,fx),\]
is called the \emph{cograph} of $f$. Obviously, for every $[0,1]_*$-set $(X,\al)$, the identity function $1_X$ is monotone, and
\[\al=(1_X)_{\nat}=1_X^{\nat}.\]
Hence, in order to simplify the notation, we abbreviate a $[0,1]_*$-set $(X,\al)$ to $X$, and write $1_X^{\nat}(x,y)$ instead of $\al(x,y)$ if no confusion arises. 

For each $q\in[0,1]$, we have a one-element $[0,1]_*$-set $\{q\}$ with 
\[1_{\{q\}}^{\nat}(q,q)=q.\]


\begin{lem} \label{phi-p-q}
For $p,q\in[0,1]$, the following statements are equivalent:
\begin{enumerate}[label={\rm(\roman*)}]
\item \label{phi-p-q:phi} There exists a morphism $\phi\colon\{p\}\oto\{q\}$ of one-element $[0,1]_*$-sets.
\item \label{phi-p-q:p} Either $p=q$, or else $p<q$ and $p$ is idempotent.
\end{enumerate}
In this case, it necessarily holds that $\phi(p,q)=p$. Therefore, a morphism between one-element $[0,1]_*$-sets, when it exists, will be denoted by
\[\overline{p}\colon\{p\}\oto\{q\}.\]
\end{lem}

\begin{proof}
\ref{phi-p-q:phi}$\implies$\ref{phi-p-q:p}: Suppose that $\phi(p,q)=r$. Then, by \ref{M1} and \ref{M4} we have
\begin{equation} \label{phi-p-q:M134}
r\leq p\wedge q\quad\text{and}\quad p\leq r*(q\ra r).
\end{equation}
Thus
\[p\leq r*(q\ra r)\leq r\leq p,\] 
and consequently $r=p$. It follows that
\[p=r\leq p\wedge q\leq q,\]
and 
\[p\leq r*(q\ra r)=p*(q\ra p)\leq p,\]
where the latter inequality forces 
\begin{equation} \label{p=p*(q-ra-p)}
p=p*(q\ra p).
\end{equation}

If $p<q$, by Lemma \ref{t-norm-rep} we have three cases:
\begin{enumerate}[label=(\alph*)]
\item \label{phi-p-q:a} There exists $s\in(p,q)$ such that $s$ is idempotent. Then $p$ is also idempotent, because
\[p=p*(q\ra p)\leq p*(s\ra p)=p*p\leq p,\]
where $s\ra p=p$ follows from \eqref{p-idem-q}.
\item \label{phi-p-q:b} There exists $[a,b]\subseteq[0,1]$ such that $a\leq p<q\leq b$ and that there exists an isomorphism $f\colon[a,b]_*\to[0,1]_{\times}$. By \eqref{p=p*(q-ra-p)} and Example \ref{MLP-def}\ref{MLP-def:P} we have 
\[f(p)=f(p)\cdot\dfrac{f(p)}{f(q)},\]
and consequently $f(p)=0$ or $f(p)=f(q)$. Since $p<q$, we conclude that $f(p)=0$. Thus $p=a$ is idempotent.
\item \label{phi-p-q:c} There exists $[a,b]\subseteq[0,1]$ such that $a\leq p<q\leq b$ and that there exists an isomorphism $g\colon[a,b]_*\to[0,1]_{*_{\L}}$. By \eqref{p=p*(q-ra-p)} and Example \ref{MLP-def}\ref{MLP-def:L} we have 
\[g(p)=g(p)*_{\L}(1-g(q)+g(p))=0\vee(2g(p)-g(q)),\]
and consequently $g(p)=0$ or $g(p)=g(q)$. Since $p<q$, we conclude that $g(p)=0$. Thus $p=a$ is idempotent.
\end{enumerate}

\ref{phi-p-q:p}$\implies$\ref{phi-p-q:phi}: In this case, it is easy to see that $\phi(p,q)=p$ satisfies \ref{M1}--\ref{M4}:
\begin{itemize}
\item \ref{M1} $p\leq p\wedge q$.
\item \ref{M2} $(q\ra q)*p*(p\ra p)\leq p$.
\item \ref{M3} $p*(p\ra p)\leq q$.
\item \ref{M4} $p\leq p*(q\ra p)$.
\end{itemize}
Therefore, $\phi\colon\{p\}\oto\{q\}$ is a morphism of one-element $[0,1]_*$-sets.
\end{proof} 

\begin{rem}
An anonymous referee suggests an alternative way to proving ``\ref{phi-p-q:phi}$\implies$\ref{phi-p-q:p}'' of Lemma \ref{phi-p-q}. In the paragraph below \eqref{p=p*(q-ra-p)}, if $p<q$, we may consider the decreasing sequence
\[\{(q\ra p)^n\}_{n\in\bbN}=\{\underbrace{(q\ra p)*\dots *(q\ra p)}\limits_{n\ \text{times}}\}_{n\in\bbN}\]
in $[0,1]$, which necessarily converges. Suppose that $s=\lim\limits_{n\rightarrow\infty}(q\ra p)^n$. Then $s$ is idempotent, since the continuity of $*$ implies that
\[s*s=\lim\limits_{n\rightarrow\infty}(q\ra p)^n*\lim\limits_{m\rightarrow\infty}(q\ra p)^m=\lim\limits_{n\rightarrow\infty}\lim\limits_{m\rightarrow\infty}(q\ra p)^{n+m}=\lim\limits_{n\rightarrow\infty}s=s.\]
Moreover, using \eqref{p=p*(q-ra-p)} and induction we deduce that 
\begin{equation} \label{p=p*(q-p)^n}
p=p*(q\ra p)^n
\end{equation}
for all $n\in\bbN$. Thus,
\begin{align*}
q\wedge s&=q*s & (s\ \text{is idempotent})\\
&=q*\lim\limits_{n\rightarrow\infty}(q\ra p)^n\\
&=\lim\limits_{n\rightarrow\infty}q*(q\ra p)^n & (*\ \text{is continuous})\\
&=\lim\limits_{n\rightarrow\infty}p*(q\ra p)^{n-1} & (\text{by Proposition \ref{divisible}\ref{divisible:u-q}})\\
&=\lim\limits_{n\rightarrow\infty}p &(\text{by \eqref{p=p*(q-p)^n}})\\
&=p. 
\end{align*}
Since $p<q$, we conclude that $p=s$, and consequently $p$ is idempotent.
\end{rem}

\section{Cauchy complete real-valued sets}

A \emph{singleton} on a $[0,1]_*$-set $X$ is a morphism $\lam\colon\{p\}\oto X$ whose domain is a one-element $[0,1]_*$-set. It is easy to see that $1_X^{\nat}(x,-)$ is a singleton for all $x\in X$. 

A $[0,1]_*$-set $X$ is said to be \emph{Cauchy complete} if it satisfies one of the equivalent conditions in the following well-known proposition:

\begin{prop} \label{Cauchy-complete-def} (See \cite[Proposition 7.1]{Stubbe2005}.)
For a $[0,1]_*$-set $X$, the following statements are equivalent:
\begin{enumerate}[label={\rm(\roman*)}]
\item \label{Cauchy-complete-def:F} Each morphism $\phi\colon A\oto X$ is the graph of a monotone function $f\colon A\to X$.
\item \label{Cauchy-complete-def:s} For each singleton $\lam\colon\{p\}\oto X$, there exists $x\in X$ such that $\lam=1_X^{\nat}(x,-)$.
\end{enumerate}
\end{prop}

In particular, in Proposition \ref{Cauchy-complete-def}\ref{Cauchy-complete-def:s} it necessarily holds that
\[1_X^{\nat}(x,x)=p,\]
because
\[p\leq\bigvee\limits_{x'\in X} 1_X^{\nat}(x',x)*(1_X^{\nat}(x',x')\ra 1_X^{\nat}(x,x'))\leq 1_X^{\nat}(x,x)\leq p,\]
where the three inequalities follow from \ref{M4}, \ref{M3} and \ref{M1}, respectively.

Let $X$ be a Cauchy complete $[0,1]_*$-set. For each $x\in X$ and each morphism $\overline{p}\colon\{p\}\oto\{1_X^{\nat}(x,x)\}$ between one-element $[0,1]_*$-sets (see Lemma \ref{phi-p-q}), 
\[1_X^{\nat}(x,-)\circ\overline{p}\colon\{p\}\oto\{1_X^{\nat}(x,x)\}\oto X\]
is also a singleton. Thus, by the Cauchy completeness of $X$, there exists $x\bullet p\in X$ such that
\begin{equation} \label{x-cdot-u-def}
1_X^{\nat}(x\bullet p,x\bullet p)=p\quad\text{and}\quad 1_X^{\nat}(x\bullet p,-)= 1_X^{\nat}(x,-)\circ\overline{p},
\end{equation}
which is necessarily unique if $X$ is separated. In fact, when $X$ is a separated Cauchy complete $[0,1]_*$-set, \eqref{x-cdot-u-def} indicates that there are bijections (cf. \eqref{X_q-def})
\begin{equation} \label{X_q-cong-RSet(q,X)}
X_q\cong\RSet(\{q\}, X)
\end{equation}
natural in $q$; that is, we have a bijection (cf. \eqref{FYd-def} below)
\begin{equation} \label{sY-def}
\sY\colon X_q\to\RSet(\{q\}, X),\quad \sY x=1_X^{\nat}(x,-)
\end{equation}
for every $q\in[0,1]$, such that for each morphism $\overline{p}\colon\{p\}\oto\{q\}$ between one-element $[0,1]_*$-sets, the square
\[\bfig
\square<1000,400>[X_q`\RSet(\{q\}, X)`X_p`\RSet(\{p\}, X);\sY`-\bullet p`-\circ\overline{p}`\sY]
\efig\]
is commutative.

\begin{prop} \label{F(x-cdot-u)}
Let $f\colon X\to Y$ be a monotone function between Cauchy complete $[0,1]_*$-sets. Then
\begin{equation} \label{f-keep-cdot}
f(x\bullet p)\cong fx\bullet p
\end{equation}
for all $x\in X$ and morphisms $\overline{p}\colon\{p\}\oto\{1_X^{\nat}(x,x)\}$, where the isomorphism is defined by \eqref{x-cong-y-def}.
\end{prop}

\begin{proof}
First, by \eqref{x-cdot-u-def} we obtain that 
\[1_Y^{\nat}(fx\bullet p,fx\bullet p)=p=1_X^{\nat}(x\bullet p,x\bullet p)=1_Y^{\nat}(f(x\bullet p),f(x\bullet p)).\]
Consequently, using \eqref{x-cdot-u-def}, \eqref{monotone-function-def} and \ref{S1}, we have
\[p=1_X^{\nat}(x\bullet p,x\bullet p)= 1_X^{\nat}(x,x\bullet p)\circ\overline{p}\leq 1_Y^{\nat}(fx,f(x\bullet p))\circ\overline{p}= 1_Y^{\nat}(fx\bullet p,f(x\bullet p))\leq 1_Y^{\nat}(fx\bullet p,fx\bullet p)=p.\]
Thus
\[1_Y^{\nat}(fx\bullet p,fx\bullet p)=1_Y^{\nat}(f(x\bullet p),f(x\bullet p))=1_Y^{\nat}(fx\bullet p,f(x\bullet p))=p,\]
and therefore, \eqref{f-keep-cdot} follows from \eqref{x-cong-y-def}.
\end{proof}

For each $[0,1]_*$-set $X$, let 
\[\CdX\coloneqq\{\lam\in\RSet(\{q\},X)\mid q\in[0,1]\}\] 
denote the set of all singletons on $X$. Then $\CdX$ is also a $[0,1]_*$-set, with
\begin{equation} \label{CdX-hom}
1_{\CdX}^{\nat}(\lam,\lam')=\bv_{x\in X}\lam(x)*(1_X^{\nat}(x,x)\ra\lam'(x))=\bv_{x\in X}\lam'(x)*(1_X^{\nat}(x,x)\ra\lam(x)),
\end{equation}
where the last equality follows from \ref{M1} and Proposition \ref{divisible}\ref{divisible:u-v-q}. Indeed, $\CdX$ is separated and Cauchy complete (see \cite[Proposition 7.12]{Stubbe2005}), called the \emph{Cauchy completion} of $X$. It is straightforward to check that for any $\lam\colon\{q\}\oto X$ in $\CdX$ and $\overline{p}\colon\{p\}\oto\{q\}$, we have
\begin{equation} \label{CdX-type}
1_{\CdX}^{\nat}(\lam,\lam)=q
\end{equation}
and
\begin{equation} \label{CdX-cdot}
\lam\bullet p=\lam\circ\overline{p}.
\end{equation}

For every $[0,1]_*$-set $X$, there is a monotone function
\begin{equation} \label{FYd-def}
\Fyd_X\colon X\to\CdX,\quad\Fyd_X x= 1_X^{\nat}(x,-).
\end{equation}
If $X$ is separated and Cauchy complete, then $X$ is isomorphic to $\CdX$ in $\RSetM$, with the isomorphism given by $\Fyd_X$, whose inverse is given by
\begin{equation} \label{Fr-def}
\Fr_X\colon\CdX\to X\quad\text{with}\quad\lam= 1_X^{\nat}(\Fr_X\lam,-).
\end{equation}

\begin{exmp} \label{Cdq-exmp}
By Lemma \ref{phi-p-q}, the Cauchy completion $\sCd\{q\}$ of a one-element $[0,1]_*$-set $\{q\}$ consists of morphisms $\overline{p}\colon\{p\}\oto\{q\}$ $(p\in[0,1])$; that is,
\begin{equation} \label{Cdq}
\sCd\{q\}=\{\overline{q}\}\cup\{\overline{p}\mid p<q\ \text{and}\ p\ \text{is idempotent}\}.
\end{equation}
Note that
\begin{equation} \label{Cdp-hom}
1_{\sCd\{q\}}^{\nat}(\overline{p},\overline{p'})=p*(q\ra p')=p\wedge p'
\end{equation}
for all $\overline{p},\overline{p'}\in\sCd\{q\}$, where the last equality holds because
\[p\wedge p'=p*p'=p*(1\ra p')\leq p*(q\ra p')=p'*(q\ra p)\leq p\wedge p'\]
if at least one of $p$ and $p'$ is idempotent, and it holds trivially if $p=p'=q$.
\end{exmp}

\begin{exmp} \label{Cdq-x-y-exmp}
Let $X$ be a two-element $[0,1]_*$-set, say $X=\{x,y\}$. Then the Cauchy completion of $X$ is
\[\sCd X=\{1_X^{\nat}(z,-)\circ\overline{p}\mid z\in\{x,y\},\ \overline{p}\colon \{p\}\oto\{1_X^{\nat}(z,z)\}\ \text{is a morphism},\ p\in[0,1]\}.\]
Indeed, if $\lam\colon\{p\}\oto X$ is a singleton, then by \ref{M4} and \ref{M1},
\[p\leq\left(\lam(x)*(1_X^{\nat}(x,x)\ra\lam(x))\right)\vee\left(\lam(y)*(1_X^{\nat}(y,y)\ra\lam(y))\right)\leq p.\]
Thus $p=\lam(x)*(1_X^{\nat}(x,x)\ra\lam(x))$ or $p=\lam(y)*(1_X^{\nat}(y,y)\ra\lam(y))$. Suppose that $p=\lam(x)*(1_X^{\nat}(x,x)\ra\lam(x))$. For $z\in\{x,y\}$, note that
\begin{align*}
\lam(z)&=(p\ra\lam(z))*p & (\text{by \ref{M1} and Proposition \ref{divisible}\ref{divisible:u-q}})\\
&=(p\ra\lam(z))*\lam(x)*(1_X^{\nat}(x,x)\ra\lam(x)) & (p=\lam(x)*(1_X^{\nat}(x,x)\ra\lam(x)))\\
&\leq 1_X^{\nat}(x,z)*(1_X^{\nat}(x,x)\ra\lam(x)) & (\text{by \ref{M3}})\\
&=(1_X^{\nat}(x,x)\ra 1_X^{\nat}(x,z))*\lam(x)) & (\text{by Proposition \ref{divisible}\ref{divisible:u-v-q}})\\
&\leq\lam(z).& (\text{by \ref{M2}})
\end{align*}
Thus
\[\lam=1_X^{\nat}(x,-)*(1_X^{\nat}(x,x)\ra\lam(x))=1_X^{\nat}(x,-)\circ\lam(x),\]
where the second equality follows from \eqref{psi-circ-phi-def}. We claim that 
\begin{equation} \label{lam-x-p-1_X-x-x}
\lam(x)\colon \{p\}\oto \{1_X^{\nat}(x,x)\}
\end{equation}
is a morphism, which would force $\lam(x)=\overline{p}$ by Lemma \ref{phi-p-q}. Indeed, \eqref{lam-x-p-1_X-x-x} satisfies \ref{M1}, \ref{M2} and \ref{M3} because $\lam$ is a singleton, and satisfies \ref{M4} because $p=\lam(x)*(1_X^{\nat}(x,x)\ra\lam(x))$, as desired.
\end{exmp}

The full subcategory of $\RSetM$ consisting of separated Cauchy complete $[0,1]_*$-sets is denoted by
\[\RCcSetM,\]
and it is a reflective subcategory of $\RSetM$ (see \cite[Proposition 7.14]{Stubbe2005}), with the reflector given by the Cauchy completion $\sCd$. The adjunction
\[\bfig
\morphism|a|/@{->}@<5pt>/<1000,0>[\RSetM\ `\ \RCcSetM;\sCd]
\morphism/@{<-^)}@<-5pt>/<1000,0>[\RSetM\ `\ \RCcSetM;]
\place(470,0)[\bot]
\efig\]
can be described by the bijections
\begin{equation} \label{RCcSetM-RSetM-ref}
\RCcSetM(\sCd X,Y)\cong\RSetM(X,Y)
\end{equation}
natural in $X\in\RSetM$ and $Y\in\RCcSetM$. In fact, the unit of this adjunction is $\Fyd$ (see \eqref{FYd-def}), and for any monotone function $f\colon X\to Y$, where $X$ is a $[0,1]_*$-set and $Y$ is a separated Cauchy complete $[0,1]_*$-set, there exists a unique monotone function (cf. \eqref{Fr-def})
\begin{equation} \label{Cf-def}
\check{f}\colon\CdX\to Y,\quad\check{f}\lam=\Fr_Y(f_{\nat}\circ\lam)
\end{equation}
such that the triangle
\[\bfig
\qtriangle<600,400>[X`\CdX`Y;\Fyd_X`f`\check{f}]
\efig\]
is commutative.


The following proposition is well known, and it allows us to explore the cartesian closedness of $\RSet$ through $\RCcSetM$ instead:

\begin{prop} \label{RSet-RCcSetM} (See \cite[Proposition 5.7(2)]{Pu2012}.)
The category $\RSet$ is equivalent to $\RCcSetM$.
\end{prop}

The adjoint equivalence between $\RSet$ and $\RCcSetM$ can be described by the bijections
\begin{equation} \label{RCcSetM-RSet-Eq}
\RCcSetM(\sCd X,Y)\cong\RSet(X,Y)
\end{equation}
natural in $X\in\RSet$ and $Y\in\RCcSetM$. Explicitly, every monotone function $f\colon\sCd X\to Y$ corresponds to a morphism 
\[(f\circ\Fyd_X)_{\nat}\colon X\oto Y\] 
of $[0,1]_*$-sets; and conversely, every morphism $\phi\colon X\oto Y$ yields to a monotone function
\begin{equation} \label{sR-def}
\sR\phi\colon\sCd X\to Y,\quad (\sR\phi)\lam=\Fr_Y(\phi\circ\lam).
\end{equation}
The naturality of \eqref{RCcSetM-RSet-Eq} in $X$ means that for each morphism $\psi\colon X'\oto X$ of $[0,1]_*$-sets, the square
\[\bfig
\square<1200,400>[\RSet(X,Y)`\RCcSetM(\sCd X,Y)`\RSet(X',Y)`\RCcSetM(\sCd X',Y);\sR`-\circ\psi`\RCcSetM(\sCd\psi,Y)`\sR]
\efig\]
is commutative, where the map $\sR$ is defined by \eqref{sR-def}. Explicitly, for every morphism $\phi\colon X\oto Y$ between $[0,1]_*$-sets,
\begin{equation} \label{hat-phi-psi=hat-phi-psi}
\sR(\phi\circ\psi)=\sR\phi\circ(\psi\circ-)\colon\sCd X'\to^{\psi\circ-}\sCd X\to^{\sR\phi}Y.
\end{equation}

\section{Cartesian closedness of the category of real-valued sets}

Recall that in a category $\CC$ with finite products, an object $Y$ is \emph{exponentiable} if the functor $-\times Y\colon\CC\to\CC$ admits a right adjoint $(-)^Y\colon\CC\to\CC$; that is, if there are bijections
\begin{equation} \label{C-X-Y-Z-cong}
\CC(X\times Y,Z)\cong\CC(X,Z^Y)
\end{equation}
natural in $X,Z\in\CC$. Moreover, $\CC$ is \emph{cartesian closed} if every $Y\in\ob\CC$ is exponentiable. Explicitly, the existence of the adjunction $(-\times Y)\dv(-)^Y$ may be characterized by its counit
\[\{\ev_Z\colon Z^Y\times Y\to Z\}_{Z\in\ob\CC},\]
called the \emph{evaluation}, with the universal property that for each morphism $f\colon X\times Y\to Z$ in $\CC$, there exists a unique morphism $\sE f\colon X\to Z^Y$ such that the triangle
\begin{equation} \label{cartesian-closed-def}
\bfig
\btriangle<700,400>[X\times Y`Z^Y\times Y`Z;\sE f\times 1_Y`f`\ev_Z]
\efig
\end{equation}
is commutative. Conversely, every morphism $g\colon X\to Z^Y$ corresponds to
\begin{equation} \label{sM-def}
\sM g=\ev_Z\circ(g\times 1_Y)\colon X\times Y\to Z.
\end{equation}
The naturality of \eqref{C-X-Y-Z-cong} in $X$ means that for each morphism $h\colon X'\to X$, the square
\[\bfig
\square<1000,400>[\CC(X,Z^Y)`\CC(X\times Y,Z)`\CC(X',Z^Y)`\CC(X'\times Y,Z);\sM`-\circ h`-\circ(h\times 1_Y)`\sM]
\efig\]
is commutative; that is, for every morphism $g\colon X\to Z^Y$,
\begin{equation} \label{f-g-times-1_Y}
\sM(g\circ h)=\sM g\circ(h\times 1_Y)\colon X'\times Y\to^{h\times 1_Y} X\times Y\to^{\sM g} Z.
\end{equation}

Note that the category $\RCcSetM$ is complete and cocomplete (see \cite[Section 7]{Pu2012}). In fact, it is straightforward to check that the product of $X,Y\in\RCcSetM$ is given by
\begin{equation} \label{X-times-Y-def}
X\times Y=\{(x,y)\mid x\in X,\ y\in Y,\ 1_X^{\nat}(x,x)=1_Y^{\nat}(y,y)\}\quad\text{and}\quad 1_{X\times Y}^{\nat}\left((x,y),(x',y')\right)=1_X^{\nat}(x,x')\wedge 1_Y^{\nat}(y,y')
\end{equation}
for all $x,x'\in X$, $y,y'\in Y$ with $1_X^{\nat}(x,x)=1_Y^{\nat}(y,y)$, $1_X^{\nat}(x',x')=1_Y^{\nat}('y,y')$.


For $Y,Z\in\RCcSetM$, if the the exponential $Z^Y$ exists in $\RCcSetM$, then there are bijections
\begin{align}
Z^Y_q&\cong\RSet(\{q\},Z^Y) & (\text{by \eqref{X_q-cong-RSet(q,X)}}) \label{Z^Y-underlying:q-Z^Y}\\
&\cong\RCcSetM(\sCd\{q\},Z^Y) & (\text{by \eqref{RCcSetM-RSet-Eq}}) \label{Z^Y-underlying:Cq-Z^Y}\\
&\cong\RCcSetM(\sCd\{q\}\times Y,Z) & (\text{by \eqref{C-X-Y-Z-cong}}) \label{Z^Y-underlying:Cq-Y-Z}
\end{align}
natural in $q$. Recall that the bijections above are given by $\sY$, $\sR$ and $\sM$, respectively (see \eqref{sY-def}, \eqref{sR-def} and \eqref{sM-def}). Thus, every $f\in Z^Y_q$ is mapped to 
\[\sM\sR\sY f\in \RCcSetM(\sCd\{q\}\times Y,Z)\]
under these bijections.

\begin{lem} \label{Z^Y-cdot}
For $Y,Z\in\RCcSetM$, suppose that the exponential $Z^Y$ exists in $\RCcSetM$. If $f\in Z^Y_q$ and $\overline{p}\colon\{p\}\oto\{q\}$ is a morphism between one-element $[0,1]_*$-sets, then
\[\sM\sR\sY(f\bullet p)=(\sM\sR\sY f)\circ((\overline{p}\circ-)\times 1_Y)\colon\sCd\{p\}\times Y\to Z.\]
\end{lem}

\begin{proof}
First, applying \eqref{hat-phi-psi=hat-phi-psi} to $\overline{p}\colon\{p\}\oto\{q\}$ and $1_{Z^Y}^{\nat}(f,-)\colon\{q\}\oto Z^Y$ yields
\begin{equation} \label{hat-Z^Y-f-p}
\sR\left(1_{Z^Y}^{\nat}(f,-)\circ\overline{p}\right)=\sR\left(1_{Z^Y}^{\nat}(f,-)\right)\circ(\overline{p}\circ -)\colon \sCd\{p\}\to Z^Y.
\end{equation}
Second, by applying \eqref{f-g-times-1_Y} to $(\overline{p}\circ -)\colon\sCd\{p\}\to\sCd\{q\}$ and $\sR\left(1_{Z^Y}^{\nat}(f,-)\right)\colon \sCd\{q\}\to Z^Y$ we obtain that
\begin{equation} \label{tilde-hat-1_Z^Y-f}
\sM\left(\sR\left(1_{Z^Y}^{\nat}(f,-)\right)\circ(\overline{p}\circ -)\right)=\sM\sR\left(1_{Z^Y}^{\nat}(f,-)\right)\circ((\overline{p}\circ -)\times 1_Y)\colon\sCd\{p\}\times Y\to Z.
\end{equation}
Thus
\begin{align*}
\sM\sR\sY(f\bullet p)&=\sM\sR\left(1_{Z^Y}^{\nat}(f\bullet p,-)\right) &(\text{by \eqref{sY-def}})\\
&=\sM\sR\left(1_{Z^Y}^{\nat}(f,-)\circ \overline{p}\right) & (\text{by \eqref{x-cdot-u-def}})\\
&=\sM\left(\sR\left(1_{Z^Y}^{\nat}(f,-)\right)\circ(\overline{p}\circ -)\right) & (\text{by \eqref{hat-Z^Y-f-p}})\\
&=\sM\sR\left(1_{Z^Y}^{\nat}(f,-)\right)\circ((\overline{p}\circ-)\times 1_Y) &(\text{by \eqref{tilde-hat-1_Z^Y-f}})\\
&=(\sM\sR\sY f)\circ((\overline{p}\circ-)\times 1_Y), &(\text{by \eqref{sY-def}})
\end{align*}
as desired.
\end{proof}


\begin{lem} \label{Z^Y-ev}
For $Y,Z\in\RCcSetM$, suppose that the exponential $Z^Y$ exists in $\RCcSetM$. Then the evaluation is given by
\[\ev_Z\colon Z^Y\times Y\to Z,\quad\ev_Z(f,y)=(\sM\sR\sY f)(\overline{q},y),\]
where $q=1_Y^{\nat}(y,y)$. 
\end{lem}

\begin{proof}
First, note that $\sM\sR\sY f\in\RCcSetM(\sCd\{q\}\times Y,Z)$ since, by \eqref{X-times-Y-def}, we have $1_{Z^Y}^{\nat}(f,f)=1_Y^{\nat}(y,y)=q$ . Thus $(\sM\sR\sY f)(\overline{q},y)$ is well defined, because $\overline{q}\in\sCd\{q\}$.

Second, note that
\[(\sR\sY f)\overline{q}=\Fr_{Z^Y}\left(1_{Z^Y}^{\nat}(f,-)\circ\overline{q}\right)=\Fr_{Z^Y}\left(1_{Z^Y}^{\nat}(f,-)\right)=\Fr_{Z^Y}\Fyd_{Z^Y}f=f,\]
where the first equality follows from \eqref{sY-def} and \eqref{sR-def}, the second equality holds because $\overline{q}\in\sCd\{q\}$ is exactly $1_{\{q\}}^{\nat}$, and the last equality follows since $Z^Y\in\RCcSetM$ indicates that $\Fr_{Z^Y}$ is the inverse of $\Fyd_{Z^Y}$ (see \eqref{FYd-def} and \eqref{Fr-def}).

Therefore, since $\sE$ and $\sM$ are inverses of each other, we conclude that
\[\ev_Z(f,y)=\ev_Z(\sR\sY f\times 1_Y)(\overline{q},y)=\ev_Z(\sE\sM\sR\sY f\times 1_Y)(\overline{q},y)=(\sM\sR\sY f)(\overline{q},y)\]
for all $(f,y)\in Z^Y\times Y$, where the last equality follows from the commutativity of the triangle \eqref{cartesian-closed-def}.
\end{proof}

For monotone functions $f\colon\sCd\{p\}\times Y\to Z$ and $g\colon\sCd\{q\}\times Y\to Z$, we define a subset $D(f,g)$ of the interval $[0,1]$, which consists of $u\in[0,p\wedge q]$ satisfying
\begin{equation} \label{Dfg-def}
\left(u*(p\ra r)*(q\ra r')\right)\wedge 1_Y^{\nat}(y,y')\leq 1_Z^{\nat}\left(f(\overline{r},y),g(\overline{r'},y')\right)
\end{equation}
for all $(\overline{r},y)\in\sCd\{p\}\times Y$ and $(\overline{r'},y')\in\sCd\{q\}\times Y$.

\begin{lem} \label{Z^Y-hom}
For $Y,Z\in\RCcSetM$, suppose that the exponential $Z^Y$ exists in $\RCcSetM$. Then
\begin{equation} \label{Z^Y-hom:f}
1_{Z^Y}^{\nat}(f,f)=p
\end{equation}
if the domain of the monotone function $\sM\sR\sY f$ is $\sCd\{p\}\times Y$, and
\begin{equation} \label{Z^Y-hom:DMRYfg}
1_{Z^Y}^{\nat}(f,g)=\bv D(\sM\sR\sY f,\sM\sR\sY g)
\end{equation}
for all $f,g\in Z^Y$.
\end{lem}

\begin{proof}
\eqref{Z^Y-hom:f} is an immediate consequence of the bijections \eqref{Z^Y-underlying:q-Z^Y}, \eqref{Z^Y-underlying:Cq-Z^Y} and \eqref{Z^Y-underlying:Cq-Y-Z}. For \eqref{Z^Y-hom:DMRYfg}, suppose that $f\in Z^Y_p$, $g\in Z^Y_q$ and $u\in[0,p\wedge q]$. Let us consider the $[0,1]_*$-set $X$ with
\begin{equation} \label{X-p-q-u}
X=\{x,x'\},\quad 1_X^{\nat}(x,x)=p,\quad 1_X^{\nat}(x',x')=q\quad\text{and}\quad 1_X^{\nat}(x,x')=u,
\end{equation}
and the map
\[h\colon X\to Z^Y\quad\text{with}\quad hx=f,\ hx'=g.\]
We claim that the following statements are equivalent, and consequently, \eqref{Z^Y-hom:DMRYfg} follows at once from \ref{Z^Y-hom:u}$\iff$\ref{Z^Y-hom:Dfg}:
\begin{enumerate}[label=(\roman*)]
\item \label{Z^Y-hom:u} $u\leq 1_{Z^Y}^{\nat}(f,g)$.
\item \label{Z^Y-hom:h} $h\colon X\to Z^Y$ is a monotone function.
\item \label{Z^Y-hom:Ch} $\check{h}\colon \sCd X\to Z^Y$ is a monotone function (defined by \eqref{Cf-def}).
\item \label{Z^Y-hom:MCh} $\sM\check{h}\colon \sCd X\times Y \to Z$ is a monotone function (defined by \eqref{Cf-def} and \eqref{sM-def}).
\item \label{Z^Y-hom:Dfg} $u\in D(\sM\sR\sY f,\sM\sR\sY g)$.
\end{enumerate}
Indeed, \ref{Z^Y-hom:u}$\iff$\ref{Z^Y-hom:h} is an immediate consequence of the definition of monotone functions (see \eqref{monotone-function-def}), while \ref{Z^Y-hom:h}$\iff$\ref{Z^Y-hom:Ch} and \ref{Z^Y-hom:Ch}$\iff$\ref{Z^Y-hom:MCh} follow from \eqref{RCcSetM-RSetM-ref} and \eqref{C-X-Y-Z-cong}, respectively. It remains to show that \ref{Z^Y-hom:MCh}$\iff$\ref{Z^Y-hom:Dfg}. 

First, by Example \ref{Cdq-x-y-exmp}, $\CdX$ consists of singletons on $X$ of the form 
\[1_X^{\nat}(x,-)\circ\overline{r}\quad\text{or}\quad 1_X^{\nat}(x',-)\circ\overline{r'},\]
where $\overline{r}\colon\{r\}\oto\{p\}$ and $\overline{r'}\colon\{r'\}\oto\{q\}$ are morphisms between one-element $[0,1]_*$-sets. Note that if $\left(1_X^{\nat}(x,-)\circ\overline{r},y\right)\in\CdX\times Y$, then $1_Y^{\nat}(y,y)=r$ (see \eqref{CdX-type} and \eqref{X-times-Y-def}). Thus
\begin{align*}
(\sM\check{h})\left(1_X^{\nat}(x,-)\circ\overline{r},y\right)&=\ev_Z(\check{h}\times 1_Y)\left(1_X^{\nat}(x,-)\circ\overline{r},y\right)&(\text{by \eqref{sM-def}})\\
&=\ev_Z\left(\check{h}\left(1_X^{\nat}(x,-)\circ\overline{r}\right),y\right)\\
&=\ev_Z\left(\Fr_{Z^Y}\left(h_{\nat}\circ 1_X^{\nat}(x,-)\circ\overline{r}\right),y\right)&(\text{by \eqref{Cf-def}})\\
&=\ev_Z\left(\Fr_{Z^Y}\left(h_{\nat}(x,-)\circ\overline{r}\right),y\right)&(\text{by \eqref{al-id}})\\ 
&=\ev_Z\left(\Fr_{Z^Y}\left(1_{Z^Y}^{\nat}(hx,-)\circ\overline{r}\right),y\right)&(\text{by \eqref{f-graph-def}})\\
&=\ev_Z\left(\Fr_{Z^Y}\left(1_{Z^Y}^{\nat}(f,-)\circ\overline{r}\right),y\right)&(hx=f)\\
&=\ev_Z\left(\Fr_{Z^Y}\left(1_{Z^Y}^{\nat}(f,-)\bullet r\right),y\right)&\left(\text{by \eqref{CdX-cdot} and}\ 1_{Z^Y}^{\nat}(f,-)\in\sCd Z^Y\right)\\
&=\ev_Z\left(\left(\Fr_{Z^Y}\left(1_{Z^Y}^{\nat}(f,-)\right)\bullet r\right),y\right)&(\text{by Proposition \ref{F(x-cdot-u)}})\\
&=\ev_Z(f\bullet r,y)&(\text{by \eqref{Fr-def}})\\
&=(\sM\sR\sY(f\bullet r))\left(\overline{1_Y^{\nat}(y,y)},y\right) & (\text{by Lemma \ref{Z^Y-ev}})\\
&=(\sM\sR\sY f)\circ((\overline{r}\circ-)\times 1_Y)\left(\overline{1_Y^{\nat}(y,y)},y\right)&(\text{by Lemma \ref{Z^Y-cdot}})\\
&=(\sM\sR\sY f)(\overline{r},y)
\end{align*}
for all $\overline{r}\colon\{r\}\oto\{p\}$ and $y\in Y$, where the last equality holds because $\overline{1_Y^{\nat}(y,y)}\colon\{r\}\oto\{r\}$ is actually the identity morphism on $\{r\}$. Similarly, we may compute that
\[(\sM\check{h})\left(1_X^{\nat}(x',-)\circ\overline{r'},y'\right)=(\sM\sR\sY g)(\overline{r'},y')\]
for all $\overline{r'}\colon\{r'\}\oto\{q\}$ and $y'\in Y$. 

Second, from \eqref{psi-circ-phi-def}, \eqref{X-p-q-u} and Proposition \ref{divisible}\ref{divisible:u-q} we see that
\begin{align*}
1_X^{\nat}(x,x)\circ\overline{r}&=1_X^{\nat}(x,x)*(p\ra r)=p*(p\ra r)=r,\\
1_X^{\nat}(x,x')\circ\overline{r}&=1_X^{\nat}(x,x')*(p\ra r)=u*(p\ra r),\\
1_X^{\nat}(x',x)\circ\overline{r'}&=1_X^{\nat}(x',x)*(q\ra r')=u*(q\ra r'),\\
1_X^{\nat}(x',x')\circ\overline{r'}&=1_X^{\nat}(x',x')*(q\ra r')=q*(q\ra r')=r'
\end{align*}
for all $\overline{r}\colon\{r\}\oto\{p\}$ and $\overline{r'}\colon\{r'\}\oto\{q\}$. Thus, it follows from \eqref{CdX-hom} that
\begin{align*}
&1_{\sCd X}^{\nat}\left(1_X^{\nat}(x,-)\circ\overline{r},1_X^{\nat}(x',-)\circ\overline{r'}\right)\\
={}&\bv_{z\in X}\left(1_X^{\nat}(x,z)\circ\overline{r}\right)*\left(1_X^{\nat}(z,z)\ra\left(1_X^{\nat}(x',z)\circ\overline{r'}\right)\right)\\
={}&\left[\left(1_X^{\nat}(x,x)\circ\overline{r}\right)*\left(1_X^{\nat}(x,x)\ra\left(1_X^{\nat}(x',x)\circ\overline{r'}\right)\right)\right]\vee\left[\left(1_X^{\nat}(x,x')\circ\overline{r}\right)*\left(1_X^{\nat}(x',x')\ra\left(1_X^{\nat}(x',x')\circ\overline{r'}\right)\right)\right]\\
={}&[r*(p\ra(u*(q\ra r')))]\vee[u*(p\ra r)*(q\ra r')]\\
={}&u*(p\ra r)*(q\ra r'),
\end{align*}
where the last equality holds because $u*(q\ra r')\leq u\leq p$ implies that
\[r*(p\ra(u*(q\ra r')))=u*(q\ra r')*(p\ra r)\]
by Proposition \ref{divisible}\ref{divisible:u-v-q}. Therefore, in conjunction with \eqref{X-times-Y-def} we conclude that the monotonicity of the function $\sM\check{h}$ is equivalent to
\[\left(u*(p\ra r)*(q\ra r')\right)\wedge 1_Y^{\nat}(y,y')\leq 1_Z^{\nat}\left((\sM\sR\sY f)(\overline{r},y),(\sM\sR\sY g)(\overline{r'},y')\right)\]
for all $(\overline{r},y)\in\sCd\{p\}\times Y$ and $(\overline{r'},y')\in\sCd\{q\}\times Y$, which establishes \ref{Z^Y-hom:MCh}$\iff$\ref{Z^Y-hom:Dfg} and thus completes the proof.
\end{proof}

\begin{exmp}
As suggested by an anonymous referee, we present here an application of Lemma \ref{Z^Y-hom} in the following simple setting. Let us consider the {\L}ukasiewicz t-norm $[0,1]_{*_{\L}}$ (see Example \ref{MLP-def}\ref{MLP-def:L}). Let $Y=\sCd\{1\}$, and assume that $Z^Y$ is the exponential of $Y,Z\in\RLCcSetM$. For $f,g\in Z^Y$ with 
\[1_{Z^Y}^{\nat}(f,f)=1\quad\text{and}\quad 1_{Z^Y}^{\nat}(g,g)=\dfrac{1}{2},\]
we may compute $1_{Z^Y}^{\nat}(f,g)$ as follows. Note that $D(\sM\sR\sY f,\sM\sR\sY g)$ consists of $u\in\Big[0,\dfrac{1}{2}\Big]$ satisfying
\begin{equation} \label{Luk-Dfg}
\left(u*(1\ra r)*\Big(\dfrac{1}{2}\ra r'\Big)\right)\wedge 1_Y^{\nat}(y,y')\leq 1_Z^{\nat}\left(f(\overline{r},y),g(\overline{r'},y')\right)
\end{equation}
for all $(\overline{r},y)\in\sCd\{1\}\times Y$ and $(\overline{r'},y')\in\sCd\Big\{\dfrac{1}{2}\Big\}\times Y$.
Since there is no non-trivial idempotent element in $[0,1]_{*_{\L}}$, we have
\[\sCd\{1\}=\left\{\overline{0}\colon\{0\}\oto\{1\},\ \overline{1}\colon\{1\}\oto \{1\}\right\}\quad\text{and}\quad\sCd\Big\{\dfrac{1}{2}\Big\}=\left\{\overline{0}\colon\{0\}\oto\Big\{\dfrac{1}{2}\Big\},\ \overline{\dfrac{1}{2}}\colon\Big\{\dfrac{1}{2}\Big\}\oto\Big\{\dfrac{1}{2}\Big\}\right\}.\]
Thus, it follows from \eqref{X-times-Y-def} that $\sCd\Big\{\dfrac{1}{2}\Big\}\times Y=\sCd\Big\{\dfrac{1}{2}\Big\}\times \sCd\{1\}$ contains only one element, i.e., 
\[\Big(\overline{0}\colon \{0\}\oto\Big\{\dfrac{1}{2}\Big\},\ \overline{0}\colon\{0\}\oto\{1\}\Big).\] 
Consequently, the left side of \eqref{Luk-Dfg} must be $0$ (because $y'$ must be $\overline{0}\colon\{0\}\oto\{1\}$, and thus $1_Y^{\nat}(y,y')=0$ by \eqref{CdX-hom}), which means that the inequality \eqref{Luk-Dfg} holds trivially. Hence $D(\sM\sR\sY f,\sM\sR\sY g)=\Big[0,\dfrac{1}{2}\Big]$, and from Lemma \ref{Z^Y-hom} we conclude that
\[1_{Z^Y}^{\nat}(f,g)=\bv D(\sM\sR\sY f,\sM\sR\sY g)=\bv\Big[0,\dfrac{1}{2}\Big]=\dfrac{1}{2}.\]
\end{exmp}

\begin{prop} \label{RCcSetM-cc}
The category $\RCcSetM$ is cartesian closed if, and only if, $*$ is the minimum t-norm on $[0,1]$.
\end{prop}

\begin{proof}
The ``if'' part is an immediate consequence of Remark \ref{OmSet} and Proposition \ref{RSet-RCcSetM}. For the ``only if'' part, we proceed by contradiction. Suppose that $*$ is not the minimum t-norm on $[0,1]$. Then, by Lemma \ref{t-norm-rep}, there exists a non-trivial closed interval 
\[[a,b]\subseteq[0,1]\] 
such that the continuous t-norm $[a,b]_*$ obtained by restricting $*$ to $[a,b]$ is either isomorphic to the product t-norm $[0,1]_{\times}$ or isomorphic to the {\L}ukasiewicz t-norm $[0,1]_{*_{\L}}$; in particular, no element in $(a,b)$ is idempotent.

Let us consider the $[0,1]_*$-set $X$ with
\[X=\{x,x'\},\quad 1_X^{\nat}(x,x)=1_X^{\nat}(x',x')=b\quad\text{and}\quad 1_X^{\nat}(x,x')=a.\]
Let $Y=\sCd\{b\}$ and $Z=\sCd X$. We show that $Z^Y$ equipped with \eqref{Z^Y-hom:DMRYfg} cannot be a $[0,1]_*$-set. 

First, using Example \ref{Cdq-exmp} and the fact that $(a,b)$ contains no idempotent element, we find that
\begin{equation} \label{Cdb-def}
\sCd\{b\}=\{\overline{b}\colon\{b\}\oto\{b\}\}\cup\{\overline{p}\colon\{p\}\oto\{b\}\mid 0\leq p\leq a\ \text{and}\ p\ \text{is idempotent}\}.
\end{equation}
Since $1_{\sCd\{b\}}^{\nat}(\overline{p},\overline{p})=p$ for all $\overline{p}\in\sCd\{b\}$ (see \eqref{CdX-type} and Lemma \ref{phi-p-q}), we see that each slice of $\sCd\{b\}$ contains only one element. Thus, by \eqref{X-times-Y-def}, $\sCd\{b\}\times Y=\sCd\{b\}\times \sCd\{b\}$ can be identified with $\sCd\{b\}$ itself. Moreover, it follows from \eqref{Cdp-hom} that
\[1_{\sCd\{b\}\times Y}^{\nat}(\overline{p},\overline{q})=p\wedge q\]
for all $p,q\in\sCd\{b\}$. Similarly, let $c=\dfrac{a+b}{2}$. Note that $\overline{c}\not\in\sCd\{b\}$, because $c\in(a,b)$ is not idempotent. Thus $\sCd\{c\}\times Y=\sCd\{c\}\times\sCd\{b\}$ can be identified with
\begin{equation} \label{Cdc-Y}
\{\overline{p}\colon\{p\}\oto\{b\}\mid 0\leq p\leq a\ \text{and}\ p\ \text{is idempotent}\},
\end{equation}
and
\[1_{\sCd\{c\}\times Y}^{\nat}(\overline{p},\overline{q})=p\wedge q\]
for all $p,q\in\sCd\{c\}$.

Second, from Example \ref{Cdq-x-y-exmp} we have
\[Z=\{1_X^{\nat}(z,-)\circ\overline{p}\mid z\in\{x,y\},\ \overline{p}\colon \{p\}\oto\{1_X^{\nat}(z,z)\}\ \text{is a morphism},\ p\in[0,b]\}.\]
Define $f,g,h\in Z^Y$ corresponding to
\begin{align*}
\sM\sR\sY f & \colon\sCd\{b\}\times Y\to Z,\quad (\sM\sR\sY f)\overline{p}=1_X^{\nat}(x,-)\circ\overline{p},\\
\sM\sR\sY g & \colon\sCd\{c\}\times Y\to Z,\quad (\sM\sR\sY g)\overline{p}=1_X^{\nat}(x,-)\circ\overline{p},\\
\sM\sR\sY h &\colon\sCd\{b\}\times Y\to Z,\quad (\sM\sR\sY h)\overline{p}=1_X^{\nat}(x',-)\circ\overline{p}.
\end{align*}
Analogously to the proof of Lemma \ref{Z^Y-hom}, we may compute that 
\begin{align*}
1_X^{\nat}(x,x)\circ\overline{p}&=1_X^{\nat}(x,x)*(b\ra p)=b*(b\ra p)=p,\\
1_X^{\nat}(x,x)\circ\overline{q}&=1_X^{\nat}(x,x)*(b\ra q)=b*(b\ra q)=q,\\
1_X^{\nat}(x,x')\circ\overline{p}&=1_X^{\nat}(x,x')*(b\ra p)=a*(b\ra p)=a*p,\\
1_X^{\nat}(x,x')\circ\overline{q}&=1_X^{\nat}(x,x')*(b\ra q)=a*(b\ra q)=a*q
\end{align*}
for all $\overline{p}\colon\{p\}\oto\{b\}$ and $\overline{q}\colon\{q\}\oto\{b\}$, where we have applied \eqref{p-idem-q} to the idempotent element $b$. Thus, since both $a$ and $b$ are idempotent,
\begin{align*}
&1_{\sCd\{b\}\times Y}^{\nat}(\overline{p},\overline{q})\\
={}&p\wedge q\\
={}&(p\wedge q)\vee(a\wedge p\wedge q)\\
={}&(p*q)\vee(a*p*a*p)\\
={}&(p*(b\ra q))\vee[a*p*(b\ra(a*q))]\\
={}&\left[\left(1_X^{\nat}(x,x)\circ\overline{p}\right)*\left(1_X^{\nat}(x,x)\ra\left(1_X^{\nat}(x,x)\circ\overline{q}\right)\right)\right]\vee\left[\left(1_X^{\nat}(x,x')\circ\overline{p}\right)*\left(1_X^{\nat}(x',x')\ra\left(1_X^{\nat}(x,x')\circ\overline{q}\right)\right)\right]\\
={}&\bv_{z\in X}\left(1_X^{\nat}(x,z)\circ\overline{p}\right)*\left(1_X^{\nat}(z,z)\ra\left(1_X^{\nat}(x,z)\circ\overline{q}\right)\right)\\
={}&1_Z^{\nat}\left(1_X^{\nat}(x,-)\circ\overline{p},1_X^{\nat}(x,-)\circ\overline{q}\right)\\
={}&1_Z^{\nat}((\sM\sR\sY f)\overline{p},(\sM\sR\sY f)\overline{q}),
\end{align*}
showing that $\sM\sR\sY f$ is a monotone function. Similarly, $\sM\sR\sY g$ and $\sM\sR\sY h$ are both monotone functions. Thus $f$, $g$ and $h$ are well defined.

For monotone functions $f\colon\sCd\{p\}\times Y\to Z$ and $g\colon\sCd\{q\}\times Y\to Z$, we define a subset $D(f,g)$ of the interval $[0,1]$, which consists of $u\in[0,p\wedge q]$ satisfying
\begin{equation}
\left(u*(p\ra r)*(q\ra r')\right)\wedge 1_Y^{\nat}(y,y')\leq 1_Z^{\nat}\left(f(\overline{r},y),g(\overline{r'},y')\right)
\end{equation}
for all $(\overline{r},y)\in\sCd\{p\}\times Y$ and $(\overline{r'},y')\in\sCd\{q\}\times Y$.

Finally, note that for $\overline{p}\in\sCd\{b\}\times Y$ and $\overline{q}\in\sCd\{c\}\times Y$, by similar calculations as above we obtain that
\[1_Z^{\nat}((\sM\sR\sY f)\overline{p},(\sM\sR\sY g)\overline{q})=1_Z^{\nat}\left(1_X^{\nat}(x,-)\circ\overline{p},1_X^{\nat}(x,-)\circ\overline{q}\right)=p\wedge q.\]
Moreover, for $\overline{p}\in\sCd\{b\}\times Y$ and $\overline{q}\in\sCd\{b\}\times Y$,
\begin{align*}
&1_Z^{\nat}((\sM\sR\sY f)\overline{p},(\sM\sR\sY h)\overline{q})\\
={}&1_Z^{\nat}\left(1_X^{\nat}(x,-)\circ\overline{p},1_X^{\nat}(x',-)\circ\overline{q}\right)\\
={}&\bv_{z\in X}\left(1_X^{\nat}(x,z)\circ\overline{p}\right)*\left(1_X^{\nat}(z,z)\ra\left(1_X^{\nat}(x',z)\circ\overline{q}\right)\right)\\
={}&\left[\left(1_X^{\nat}(x,x)\circ\overline{p}\right)*\left(1_X^{\nat}(x,x)\ra\left(1_X^{\nat}(x',x)\circ\overline{q}\right)\right)\right]\vee\left[\left(1_X^{\nat}(x,x')\circ\overline{p}\right)*\left(1_X^{\nat}(x',x')\ra\left(1_X^{\nat}(x',x')\circ\overline{q}\right)\right)\right]\\
={}&[p*(b\ra(a*q))]\vee[a*p*(b\ra q)]\\
={}&(p*a*q)\vee(a*p*q)\\
={}&a*p*q\\
={}&a\wedge p\wedge q;
\end{align*}
and similarly, for $\overline{p}\in\sCd\{c\}\times Y$ and $\overline{q}\in\sCd\{b\}\times Y$,
\[1_Z^{\nat}((\sM\sR\sY g)\overline{p},(\sM\sR\sY h)\overline{q})=1_Z^{\nat}\left(1_X^{\nat}(x,-)\circ\overline{p},1_X^{\nat}(x',-)\circ\overline{q}\right)=a\wedge p\wedge q.\]
Therefore:
\begin{itemize}
\item $D(\sM\sR\sY f,\sM\sR\sY g)$ consists of those $u\in[0,c]$ satisfying
\begin{equation} \label{D(MRYf,MRYg)}
(u*(b\ra p)*(c\ra q))\wedge(p\wedge q)\leq p\wedge q
\end{equation}
for all $\overline{p}\in\sCd\{b\}\times Y$ and $\overline{q}\in\sCd\{c\}\times Y$. Since \eqref{D(MRYf,MRYg)} always holds, we have $D(\sM\sR\sY f,\sM\sR\sY g)=[0,c]$.
\item $D(\sM\sR\sY g,\sM\sR\sY h)$ consists of those $u\in[0,c]$ satisfying
\begin{equation} \label{D(MRYg,MRYh)}
(u*(c\ra p)*(b\ra q))\wedge (p\wedge q)\leq a\wedge p\wedge q
\end{equation}
for all $\overline{p}\in\sCd\{c\}\times Y$ and $\overline{q}\in\sCd\{b\}\times Y$. Since $p\leq a$ by \eqref{Cdc-Y}, the inequality \eqref{D(MRYg,MRYh)} also holds trivially. Thus $D(\sM\sR\sY g,\sM\sR\sY h)=[0,c]$.
\item $D(\sM\sR\sY f,\sM\sR\sY h)$ consists of those $u\in[0,b]$ satisfying 
\begin{equation} \label{D(MRYf,MRYh)}
(u*(b\ra p)*(b\ra q))\wedge(p\wedge q)\leq a\wedge p\wedge q
\end{equation}
for all $\overline{p},\overline{q}\in\sCd\{b\}\times Y$. Note that 
\[(u*(b\ra p)*(b\ra q))\wedge(p\wedge q)=(u*p*q)\wedge p\wedge q=u\wedge p\wedge q,\]
where the last equality follows from the fact that $p$ and $q$ are both idempotent (see \eqref{Cdb-def}, and note that $b$ is also idempotent). Thus, the inequality \eqref{D(MRYf,MRYh)} becomes
\[u\wedge p\wedge q\leq a\wedge p\wedge q,\]
which holds if, and only if, $u\leq a$ (for the ``only if'' part, just note that the arbitrariness of $\overline{p},\overline{q}\in\sCd\{b\}\times Y$ allows us to choose $p=q=b$). It follows that $D(\sM\sR\sY f,\sM\sR\sY h)=[0,a]$.
\end{itemize}
Consequently, by Lemma \ref{Z^Y-hom} we obtain that
\[1_{Z^Y}^{\nat}(g,g)=c,\quad 1_{Z^Y}^{\nat}(f,g)=1_{Z^Y}^{\nat}(g,h)=\bv[0,c]=c\quad\text{and}\quad 1_{Z^Y}^{\nat}(f,h)=\bv[0,a]=a.\]
Hence
\[1_{Z^Y}^{\nat}(g,h)*\left(1_{Z^Y}^{\nat}(g,g)\ra 1_{Z^Y}^{\nat}(f,g)\right)=c*(c\ra c)=c>a= 1_{Z^Y}^{\nat}(f,h),\]
which contradicts to \ref{S3}.
\end{proof}

Therefore, the main result of this paper arises as an immediate consequence of Propositions \ref{RSet-RCcSetM} and \ref{RCcSetM-cc}:

\begin{thm} \label{QSet-topos-frame}
The category $\RSet$ is cartesian closed if, and only if, $*$ is the minimum t-norm on $[0,1]$.
\end{thm}

In particular, as an immediate corollary, we reproduce the result of \cite{Hu2025} in the framework of real-valued sets:

\begin{cor} (See \cite[Theorem 5.6]{Hu2025}.)
The category $\RSet$ is a topos if, and only if, $*$ is the minimum t-norm on $[0,1]$.
\end{cor}

\section*{Acknowledgement}

The authors acknowledge the support of National Natural Science Foundation of China (No. 12071319) and the Fundamental Research Funds for the Central Universities (No. 2021SCUNL202). The authors are grateful for helpful remarks received from Professor Hongliang Lai, Professor Dexue Zhang and the anonymous referees.






\end{document}